\documentclass[12pt]{article}

\usepackage{amsmath,amssymb,amsthm}

\advance \topmargin by -\headheight
\advance \topmargin by -\headsep
\evensidemargin \oddsidemargin
\author{J.-P. Allouche \\
CNRS, LRI, B\^atiment 490 \\
F-91405 Orsay Cedex (France)\\
{\tt allouche@lri.fr} \\
\and
J. Sondow \\
209 West 97th Street \\
New York, NY10025 (USA) \\
{\tt jsondow@alumni.princeton.edu} \\
}

\title{Infinite products with strongly $B$-multiplicative exponents}

\date{ }

\def \proof{\bigbreak\noindent{\it Proof.\ \ }}

\def \endpf{{\ \ \ $\blacksquare$}}
\def \nodiv{{|\kern-4.2pt/}}

\newtheorem{theorem}{Theorem}
\newtheorem{lemma}{Lemma}
\newtheorem{corollary}{Corollary}
\newtheorem{proposition}{Proposition}
\newtheorem*{addendum}{Addendum}

\theoremstyle{definition}

\newtheorem{remark}{Remark}
\newtheorem{example}{Example}
\newtheorem{definition}{Definition}

\begin{document}

\maketitle

\begin{center}
{\it To Professor K\'atai on the occasion of his 70th birthday}
\end{center}

\medskip

\begin{abstract}
Let $N_{1,B}(n)$ denote the number of ones in the $B$-ary expansion 
of an integer $n$. Woods introduced the infinite product
$P :=\prod_{n \geq 0} \left(\frac{2n+1}{2n+2}\right)^{(-1)^{N_{1,2}(n)}}$
and Robbins proved that $P = 1/\sqrt{2}$. Related products were studied by
several authors. We show that a trick for proving that $P^2 = 1/2$ (knowing
that $P$ converges) can be extended to evaluating new products with (generalized)
strongly $B$-multiplicative exponents. A simple example is
$$
\prod_{n \geq 0} \left(\frac{Bn+1}{Bn+2}\right)^{(-1)^{N_{1,B}(n)}} =
\frac{1}{\sqrt B}.
$$

\medskip

\noindent
{\it MSC}:  11A63, 11Y60.
\end{abstract}

\section{Introduction}

In 1985 the following infinite product, for which no closed expression is known,
appeared in \cite[p.~193 and p.~209]{FlaMar}:
$$
R := \prod_{n \geq 1} \left(\frac{(4n+1)(4n+2)}{4n(4n+3)}\right)^{\varepsilon(n)}
$$
where $(\varepsilon(n))_{n \geq 0}$ is the $\pm 1$ Prouhet-Thue-Morse sequence,
defined by 
$$
\varepsilon(n) = (-1)^{N_{1,2}(n)}
$$
with $N_{1,2}(n)$ being the number of ones in the binary expansion of $n$. 
(For more on the Prouhet-Thue-Morse sequence, see for example \cite{AS99}.)

On the one hand, it is not difficult to see that $R=\frac{3}{2Q}$, where
$$
Q := \prod_{n \geq 1} \left(\frac{2n}{2n+1}\right)^{\varepsilon(n)}.
$$
Namely, splitting the simpler product into even and odd indices 
and using the relations $\varepsilon(2n) = \varepsilon(n)$ and 
$\varepsilon(2n+1) = -\varepsilon(n)$, we get
$$
Q = \left(\prod_{n \geq 1} \left(\frac{4n}{4n+1}\right)^{\varepsilon(n)}\right)
\left(\prod_{n \geq 0} \left(\frac{4n+2}{4n+3}\right)^{-\varepsilon(n)}\right)
= \frac{3}{2} 
\prod_{n \geq 1} \left(\frac{4n(4n+3)}{(4n+1)(4n+2)}\right)^{\varepsilon(n)}
= \frac{3}{2R}.
$$
(Note that, whereas the logarithm of $R$ is an absolutely convergent series, 
the logarithm of $Q$ -- and similarly the logarithm of the product $P$ below 
-- is a conditionally convergent series, as can be seen by partial summation, 
using the fact that the sums $\sum_{0 \leq k \leq n} \varepsilon(k)$ only take 
the values $+1$, $0$ and $-1$, hence are bounded.)

On the other hand, the product $Q$ reminds us of the Woods-Robbins
product \cite{Woods, Robbins}
$$
P := \prod_{n \geq 0} \left(\frac{2n+1}{2n+2}\right)^{\varepsilon(n)} 
= \frac{1}{\sqrt{2}}
$$
(generalized for example in \cite{Shallit, AllCoh, ACMS, AMFP, AS89, Sondow}). 

In 1987 during a stay at the University of Chicago, the first author, convinced 
that the computation of the infinite product $Q$ should not resist the even-odd
splitting techniques he was using with J.~Shallit, discovered the following trick. 
First write $QP$ as
$$
QP =   \left(\frac{1}{2}\right)^{\varepsilon(0)} 
\prod_{n \geq 1}\left(\frac {2n}{2n+1} \cdot \frac {2n+1}{2n+2}\right)^{\varepsilon(n)}
= \frac{1}{2} \prod_{n \geq 1} \left(\frac{n}{n+1}\right)^{\varepsilon(n)}.
$$
Now split the indices as we did above, obtaining
$$
\prod_{n \geq 1} \left(\frac{n}{n+1}\right)^{\varepsilon(n)}
=
\left(\prod_{n \geq 1} \left(\frac{2n}{2n+1}\right)^{\varepsilon(n)}\right)
\left(\prod_{n \geq 0} \left(\frac{2n+1}{2n+2}\right)^{-\varepsilon(n)}\right)
= QP^{-1}.
$$
This gives $QP = \frac{1}{2Â} Q P^{-1}$: as the hope of computing $Q$
fades, the trick at least yields an easy way to compute $P = 1/\sqrt 2$.
By extending this trick to $B$-ary expansions, the second author \cite{Sondow}
found the generalization of $P = 1/\sqrt{2}$ given in Corollary~\ref{sumofdigits} 
of Section~\ref{section-sum}.

It happens that the sequence $(\varepsilon(n))_{n \geq 0}$ is strongly
$2$-multiplicative (see Definition~\ref{defmult} in the next section). 
The purpose of this paper is to extend the trick to products with more 
general exponents. For example, we prove the following.

{\it Let $B > 1$ be an integer. For $k = 1, \ldots, B-1$ define $N_{k,B}(n)$ 
to be the number of occurrences of the digit $k$ in the $B$-ary expansion of 
the integer $n$. Also, let 
$$ 
s_B(n) := \sum_{0 < k <B} kN_{k,B}(n)
$$
be the sum of the $B$-ary digits of $n$, and let $q > 1$ be an integer. Then
$$
\prod_{n \geq 0} \left(\frac{Bn+k}{Bn+k+1}\right)^{(-1)^{N_{k,B}(n)}}
= \frac{1}{\sqrt{B}},
$$
$$
\prod_{n \geq 0}
\prod_{\genfrac{}{}{0cm}{2}{0 < k < B}{k \not \equiv 0 \bmod q}}
\left(\frac{Bn+k}{Bn+k+1}\right)^{\sin \tfrac{\pi k}{q}
\sin \tfrac{\pi(2s_B(n)+k)}{q}} = \frac{1}{\sqrt{B}}, 
$$
and
$$
\prod_{n \geq 0}
\prod_{\genfrac{}{}{0cm}{2}{0 < k < B}{k \not \equiv 0 \bmod q}}
\left(\frac{Bn+k}{Bn+k+1}\right)^{\sin \tfrac{\pi k}{q}
\cos \tfrac{\pi(2s_B(n)+k)}{q}} = 1.
$$
}

\medskip

Note that the use of the trick is not necessarily the only way to compute 
products of this type: real analysis is used for computing $P$ in 
\cite{Robbins} and to compute products more general than $P$ in 
\cite{Shallit}; the core of \cite{AllCoh} is the use of Dirichlet series, 
while \cite{ACMS} deals with complex power series and the second part of 
\cite{AMFP} with real integrals. It may even happen that, in some cases, 
the use of the trick gives less general results than other methods. For 
example, in Remark~\ref{trickson} we show that Corollary~\ref{sumofdigits} 
or \cite{Sondow} can also be obtained as an easy consequence 
of \cite[Theorem~1]{ACMS}.

\section{Strongly $B$-multiplicative sequences}\label{infprod}

We recall the classical definition of a strongly $B$-multiplicative sequence. 
(For this and for the definitions of $B$-multiplicative, $B$-additive, and 
strongly $B$-additive, see \cite{BelSha, Gelfond, Delange, Mendes, Grabner}.) 

\begin{definition}\label{defmult}
Let $B \geq 2$ be an integer.  A sequence of complex numbers 
$(u(n))_{n \geq 0}$ is {\em strongly $B$-multiplicative} if $u(0)=1$ 
and, for all $n \geq 0$ and all $k \in \{0, 1, \ldots, B-1\}$,
$$
u(Bn + k) = u(n) u(k).
$$
\end{definition}

\begin{example}
If $z$ is any complex number, then the sequence $u$ defined by $u(0) :=1$ 
and $u(n) := z^{s_B(n)}$ for $n \geq 1$ is strongly $B$-multiplicative.
\end{example}

\begin{remark}\label{remmult}
If we do not impose the condition $u(0)=1$ in Definition~\ref{defmult}, then 
either $u(0)=1$ holds, or the sequence $(u(n))_{n \geq 0}$ must be identically $0$. 
To see this, note that the relation $u(Bn+k) = u(n)u(k)$ implies, with $n=k=0$, 
that $u(0) = u(0)^2$. Hence $u(0)=1$ or $u(0)=0$. If $u(0)=0$, then taking $n=0$ 
in the relation gives $u(k)=0$ for all $k \in \{0, 1, \ldots, B-1\}$, which by (1) 
implies $u(n)=0$ for all $n \geq 0$.
\end{remark}

\begin{proposition}\label{remmult-prop}
If the sequence $(u(n))_{n \geq 0}$ is strongly $B$-multiplicative,
and if the $B$-ary expansion of $n \geq 1$ is $n = \sum_j e_j(n) B^j$,
then $u(n) = \prod_j u(e_j(n))$. In particular, the only strongly 
$B$-multiplicative sequence with $u(1) = u(2) = \cdots = u(B-1) = \theta$,
where $\theta = 0$ or $1$, is the sequence $1, \theta, \theta, \theta, \ldots$.
\end{proposition}

\proof Use induction on the number of base $B$ digits of $n$. \endpf 

\bigskip

We now generalize the notion of a strongly $B$-multiplicative 
sequence different from $1, 0, 0, 0, \ldots$. 

\begin{definition}~\label{H}
Let $B \geq 2$ be an integer. A sequence of complex numbers $(u(n))_{n \geq 0}$ 
{\it satisfies Hypothesis ${\cal H}_B$} if there exist an integer $n_0 \geq B$ 
and complex numbers $v(0), v(1), \ldots, v(B-1)$ such that $u(n_0) \neq 0$ and,
for all $n \geq 1$ and all $k = 0, 1, \ldots, B-1$,
$$
u(Bn + k) = u(n) v(k).
$$
\end{definition}

\begin{proposition}\label{elementary}

\ { }

{\rm (1)} If a sequence $(u(n))_{n \geq 0}$ satisfies Hypothesis ${\cal H}_B$,
then the values $v(0), v(1), \ldots, v(B-1)$ are uniquely determined.

{\rm (2)} A sequence $(u(n))_{n \geq 0}$ has $u(0)=1$ and satisfies Hypothesis 
${\cal H}_B$ with $u(Bn+k) = u(n)v(k)$ not only for $n \geq 1$ but also 
for $n=0$, if and only if the sequence is strongly $B$-multiplicative and
not equal to $1, 0, 0, 0, \ldots$. In that case, $v(k) = u(k)$ for 
$k = 0,1, \ldots, B-1$.

\end{proposition}

\proof If the sequence $(u(n))_{n \geq 0}$ satisfies Hypothesis~${\cal H}_B$, 
then $v(k) = u(Bn_0+k)/u(n_0)$ for $k = 0, 1, \ldots, B-1$. This implies (1).

To prove the ``only if'' part of (2), take $n=0$ in the relation 
$u(Bn+k) = u(n)v(k)$, yielding $u(k) = u(0) v(k) = v(k)$ for 
$k = 0, 1, \ldots, B-1$. Hence $u(Bn+k) = u(n)u(k)$ for all 
$n \geq 0$ and $k = 0, 1, \ldots, B-1$. 
Thus $(u(n))_{n \geq 0}$ is strongly $B$-multiplicative. Since 
$u(n_0) \neq 0$ for some $n_0 \geq B$, the sequence is not 
$1, 0, 0, 0, \ldots$

Conversely, suppose that $(u(n))_{n \geq 0}$ is strongly $B$-multiplicative 
and is not $1, 0, 0, 0, \ldots$ Then there exists an integer
$\ell_0 \geq 1$ such that $u(\ell_0) \neq 0$. Hence $n_0 := B\ell_0 \geq B$ and
$u(n_0) = u(B\ell_0) = u(\ell_0)u(0) = u(\ell_0) \neq 0$. Defining $v(k) := u(k)$
for $k = 0, 1, \ldots, B-1$, we see that $(u(n))_{n \geq 0}$
satisfies Hypothesis ${\cal H}_B$, and the proposition follows.  \endpf

\begin{example}\label{special}
We construct a sequence which satisfies Hypothesis ${\cal H}_B$ but is 
not strongly $B$-multipli\-ca\-ti\-ve. Let $z$ be a complex number, with 
$z \notin \{0,1\}$, and define $u(n) := z^{N_{0,B}(n)}$, where 
$N_{0,B}(n)$ counts the number of zeros in the $B$-ary expansion of $n$ 
for $n > 0$, and $N_{0,B}(0) := 0$ (which corresponds to representing zero 
by the empty sum, that is, the empty word). Note that for all $n \geq 1$ 
the relation $N_{0,B}(Bn) = N_{0,B}(n)+1 $ holds, and for all 
$k \in \{1, 2, \ldots, B-1\}$ and all $n \geq 0$ the relation 
$N_{0,B}(Bn+k) = N_{0,B}(n) = N_{0,B}(n) + N_{0,B}(k)$ holds. Hence the nonzero
sequence $(u(n))_{n \geq 0}$ satisfies Hypothesis ${\cal H}_B$, with 
$v(0) := z$ and $v(k) := 1 = u(k)$ for $k = 1, 2, \ldots, B-1$. But the 
sequence is not strongly $B$-multiplicative: $u(B \times 1 + 0) = z
\neq 1 = u(1) u(0)$. 
\end{example}

\begin{remark}\label{modifiedcountofzeroes}
The alternative definition $N_{0,B}(0) := 1$ (which would correspond to representing 
zero by the single digit $0$ instead of by the empty word) would also not lead to a 
strongly $B$-multiplicative sequence $u$, since then $u(0) = z \neq 1$, which 
does not agree with Definition~\ref{defmult} (see also Remark~\ref{remmult}). On the 
other hand, the new sequence would still satisfy Hypothesis ${\cal H}_B$, with the 
same values $v(k)$, as the same proof shows, since $u(0)$ does not appear in it.
\end{remark}

\section{Convergence of infinite products}\label{convergence}

Inspired by the Woods-Robbins product $P$, we want to study products 
given in the following lemma.

\begin{lemma}\label{conv}
Let $B > 1$ be an integer.
Let $(u(n))_{n \geq 0}$ be a sequence of complex numbers with $|u(n)| \leq 1$
for all $n \geq 0$. Suppose that it satisfies Hypothesis ${\cal H}_B$ with
$|v(k)| \leq 1$ for all $k \in \{0, 1, \ldots, B-1\}$, and that
$|\sum_{0 \leq k < B} v(k)| < B$.
Then for each $k \in \{0, 1, \ldots, B-1\}$, the infinite product
$$
\prod_{n \geq \delta_k} \left(\frac{Bn+k}{Bn+k+1}\right)^{u(n)}
$$ 
converges, where $\delta_k$ ---a special case of the Kronecker delta--- 
is defined by 
$$
\delta_k :=
\begin{cases}
             0 &\mbox{\rm if } k \neq 0 \\
             1 &\mbox{\rm if } k = 0.
\end{cases}
$$
\end{lemma}

\proof For $N = 1, 2, \ldots$, let
$$
F(N) := \sum_{0 \leq n < N} u(n).
$$
Also define for $j = 1, 2, \ldots, B-1$
$$
G(j) := \sum_{0 \leq n < j} v(n)
$$
and set $G(0) :=0$.
Then, for each $b \in \{0, 1, \ldots, B-1\}$, and for every $N \geq 1$,
$$
\begin{array}{lll}
F(BN+b) &=& \displaystyle\sum_{0 \leq n < BN} u(n) + \sum_{BN \leq n < BN+b} u(n) \\
&=& \displaystyle\sum_{0 \leq n < N} \sum_{0 \leq \ell < B} u(Bn+\ell) 
+  \sum_{0 \leq \ell < b} u(BN+\ell) \\
&=& \displaystyle\sum_{0 \leq \ell < B} u(\ell) + 
\sum_{1 \leq n < N} \sum_{0 \leq \ell < B} u(n) v(\ell) 
+ u(N) \sum_{0 \leq \ell < b} v(\ell).
\end{array}
$$
Hence, using $|u(N)| \leq 1$ and $|G(b)| \leq B-1 < B$,
$$ 
\begin{array}{lll}
|F(BN+b)| &=& |F(B) + (F(N)-u(0)) G(B) + u(N) G(b) | \\
&<& |F(B) - u(0)G(B)| + |F(N)||G(B)| + B.
\end{array}
$$
This gives the case $d=1$ of the following inequality, which holds for $d \geq 1$
and $e_t \in \{0, 1, \ldots, B-1\}$, and which is proved by induction on $d$ using
$|F(e_t)| \leq B$:
$$
\left|F\left(\sum_{0 \leq t \leq d} e_t B^t\right)\right| < 
|F(B) - u(0)G(B)|\left(1 + \sum_{1 \leq t \leq d-1} |G(B)|^t\right)
+ B \left(1 + \sum_{1 \leq t \leq d} |G(B)|^t\right).
$$
Hence
$$
\left|F\left(\sum_{0 \leq t \leq d} e_t B^t\right)\right| <
\begin{cases}
B(3d+1) &\mbox{\rm if } |G(B)| \leq 1, \\
\displaystyle 3B\frac{|G(B)|^{d+1}-1}{|G(B)|-1} &\mbox{\rm if } |G(B)| > 1.
\end{cases}
$$
This implies that for some constant $C = C(B,v)$, and for every $N$ large enough, 
$$
|F(N)| <
\begin{cases}
C \log N &\mbox{\rm if } |G(B)| \leq 1, \\
C |G(B)|^{\frac{\log N}{\log B}} = C N^{\frac{\log |G(B)|}{\log B}}
&\mbox{\rm if } |G(B)| > 1.
\end{cases}
$$
Since $|G(B)| < B$ by hypothesis, we can define $\alpha \in (0,1)$ by 
$$
\alpha := \begin{cases}
\frac{1}{2} &\mbox{\rm if } |G(B)| \leq 1, \\
\frac{\log |G(B)|}{\log B} 
&\mbox{\rm if } |G(B)| > 1.
\end{cases}
$$
Hence for every $N$ large enough $|F(N)| < C N^{\alpha}$.
It follows, using summation by parts, that the series 
$\sum_n u(n) \log \frac{Bn+k}{Bn+k+1}$ converges, hence 
the lemma.  \endpf

\begin{remark} 

\ { }

(1) Here and in what follows, expressions of the form $a^z$, where $a$
is a positive real number and $z$ a complex number, are defined by
$a^z := e^{z \log a}$, and $\log a$ is real.

(2) For more precise estimates of summatory functions of (strongly)
$B$-multiplicative sequences, see for example \cite{Delange, Grabner}.
(In \cite{Grabner} strongly $B$-multiplicative sequences are called
completely $B$-multiplicative.)
\end{remark}

\section{Evaluation of infinite products}\label{products}

This section is devoted to computing some infinite products with
exponents that satisfy Hypothesis ${\cal H}_B$, including some whose 
exponents are strongly $B$-multiplicative.

\subsection{General results}

\begin{theorem}\label{general}
Let $B > 1$ be an integer. Let $(u(n))_{n \geq 0}$ be a sequence of complex
numbers with $|u(n)| \leq 1$ for all $n \geq 0$. Suppose that $u$ satisfies 
Hypothesis ${\cal H}_B$, with complex numbers $v(0), v(1), \ldots, v(B-1)$ 
such that $|v(k)| \leq 1$ for $k \in \{0, 1, \ldots B-1\}$ and 
$|\sum_{0 \leq k < B} v(k)| < B$. Then the following relation between 
nonempty products holds:
$$
\prod_{\genfrac{}{}{0cm}{2}{0 \leq k < B}{v(k) \neq 1}} 
\prod_{n \geq \delta_k}
\left(\frac{Bn+k}{Bn+k+1}\right)^{u(n)(1-v(k))}
= \frac{1}{B^{u(0)}} \prod_{0 < k < B} 
\left(\frac{k}{k+1}\right)^{u(k)-u(0)v(k)}.
$$
\end{theorem}

\proof The condition $|\sum_{0 \leq k < B} v(k)| < B$ prevents $v$ from being 
identically equal to $1$ on $\{0, 1, \ldots, B-1\}$, so the left side of the 
equation is not empty. Since $B > 1$, so is the right.

We first show that 
\renewcommand{\theequation}{\fnsymbol{equation}}
\begin{equation}\label{split}
\prod_{0 \leq k < B} 
\prod_{n \geq \delta_k} \left(\frac{Bn+k}{Bn+k+1}\right)^{u(n)}
= \frac{1}{B^{u(0)}} \prod_{n \geq 1} \left(\frac{n}{n+1}\right)^{u(n)}
\end{equation}
(note that, by Lemma~\ref{conv}, all the products converge).
To see this, write the left side as
$$
\left(\frac{1}{2} \frac{2}{3} \cdots \frac{B-1}{B}\right)^{u(0)}\prod_{n \geq 1}
\left(\frac{Bn}{Bn+1} \frac{Bn+1}{Bn+2} \cdots \frac{Bn+B-1}{Bn+B}\right)^{u(n)}
$$
and use telescopic cancellation.
Now, splitting the product on the right side of $(\ref{split})$ according 
to the values of $n$ modulo $B$ gives
$$
\begin{array}{lll}
\displaystyle \prod_{n \geq 1} \left(\frac{n}{n+1}\right)^{u(n)} &=&
\displaystyle \prod_{0 \leq k < B} \prod_{n \geq \delta_k}
\left(\frac{Bn+k}{Bn+k+1}\right)^{u(Bn+k)} \\
&=&\displaystyle\prod_{0 < k < B} \left(\frac{k}{k+1}\right)^{u(k)}
\prod_{0 \leq k < B} \prod_{n \geq 1}
\left(\frac{Bn+k}{Bn+k+1}\right)^{u(n)v(k)} \\
&=&\displaystyle\prod_{0 < k < B} \left(\frac{k}{k+1}\right)^{u(k)-u(0)v(k)}
\prod_{0 \leq k < B} \prod_{n \geq \delta_k}
\left(\frac{Bn+k}{Bn+k+1}\right)^{u(n)v(k)}.
\end{array}
$$
Using $(\ref{split})$ and the fact that convergent infinite products 
are nonzero, the theorem follows.  \endpf

\begin{example}\label{firstcounting0}
As in Example~\ref{special}, the sequence $u$ defined by
$u(n) = z^{N_{0,B}(n)}$, with $z \notin \{0,1\}$, satisfies
Hypothesis ${\cal H}_B$, and $\sum_{0 \leq k < B} v(k) = z+B-1$. 
If furthermore $|z| \leq 1$, then 
$$
\prod_{n \geq 1} \left(\frac{Bn}{Bn+1}\right)^{(1-z)z^{N_{0,B}(n)}} = B.
$$
\end{example}

\begin{corollary}\label{generalmult}
Fix an integer $B > 1$. If $(u(n))_{n \geq 0}$ is strongly $B$-multiplicative, 
satisfies $|u(n)| \leq 1$ for all $n \geq 0$, and is not equal to either of the 
sequences $1, 0, 0, 0, \ldots$ or $1, 1, 1, \ldots$, then
$$
\prod_{n \geq 0} \prod_{\genfrac{}{}{0cm}{2}{0 < k < B}{u(k) \neq 1}}
\left(\frac{Bn+k}{Bn+k+1}\right)^{u(n)(1-u(k))} = \frac{1}{B}.
$$
\end{corollary}

\proof Using Theorem~\ref{general} and Proposition~\ref{elementary}~part (2) it 
suffices to prove that $|\sum_{0 \leq k < B} u_k| < B$. Since $|u_n| \leq 1$
for all $n \geq 0$, we have $|\sum_{0 \leq k < B} u_k| \leq B$. From the equality
case of the triangle inequality, it thus suffices to prove that the numbers
$u_0, u_1, \ldots, u_{B-1}$ are not all equal to a same complex number $z$ with $|z|=1$.
If they were, then, since $u_0 = 1$, we would have $u_0 = u_1 = \ldots = u_{B-1} = 1$. 
Hence $(u(n))_{n \geq 0} = 1, 1, 1, \ldots$ from Proposition~\ref{remmult-prop}, a
contradiction.  \endpf

\begin{addendum}
Theorem~\ref{general} and Corollary~\ref{generalmult} can be strengthened,
as follows.

{\rm (1)} If $B$, $u$, and $v$ satisfy the hypotheses of Theorem~\ref{general}, then
$$
\sum_{\genfrac{}{}{0cm}{2}{0 \leq k < B}{v(k) \neq 1}} (1-v(k))
\sum_{n \geq \delta_k} u(n) \log \frac{Bn+k}{Bn+k+1}
= - u(0) \log B + \sum_{0 < k < B}
(u(k)-u(0)v(k)) \log \frac{k}{k+1}.
$$

\medskip

{\rm (2)} If $B$ and $u$ satisfy the hypotheses of Corollary~\ref{generalmult}, then
$$
\sum_{n \geq 0} \sum_{\genfrac{}{}{0cm}{2}{0 < k < B}{u(k) \neq 1}}
u(n)(1-u(k)) \log \frac{Bn+k}{Bn+k+1} = - \log B.
$$
\end{addendum}

\proof Write the proofs of Theorem~\ref{general} and Corollary~\ref{generalmult} 
additively instead of multiplicatively.
\endpf

\begin{remark} The Addendum cannot be proved by just taking logarithms in the 
formulas in Theorem~\ref{general} and Corollary~\ref{generalmult}.
To illustrate the problem, note that while
$$
\prod_{n \geq 0} e^{\tfrac{(-1)^n 8i}{2n+1}} = 1
$$
(because the product converges to $e^{2 \pi i}$), the log equation is false:
$$
\sum_{n \geq 0} \frac{(-1)^n 8i}{2n+1} = 2 \pi i \neq  0 = \log 1.
$$
\end{remark}

\begin{example}\label{counting0}
With the same $u$ and $z$ as in Example~\ref{firstcounting0}, Addendum (1) yields
$$
\sum_{n \geq 1} z^{N_{0,B}(n)} \log \frac{Bn}{Bn+1} = \frac{\log B}{z-1}.
$$
Hence
$$
\prod_{n \geq 1} \left(\frac{Bn}{Bn+1}\right)^{z^{N_{0,B}(n)}}
= B^{\tfrac{1}{z-1}}.
$$
(Note the similarity between this product and the one in 
Example~\ref{firstcounting0}. Neither implies the other, 
but of course the preceding log equation implies both.)

If we modify the sequence $u$ as in Remark~\ref{modifiedcountofzeroes}, we get 
the same two formulas, because the value $N_{0,B}(0)$ does not appear in them.
\end{example}

\begin{corollary}\label{rootofunity}
Fix integers $B, q, p$ with $B > 1$, $q > p > 0$, and $B \equiv 1 \bmod q$. 
Then
$$
\prod_{n \geq 0}
\prod_{\genfrac{}{}{0cm}{2}{0 < k < B}{k \not \equiv 0 \bmod q}}
\left(\frac{Bn+k}{Bn+k+1}\right)^{\sin \tfrac{\pi k p}{q}
\sin \tfrac{\pi(2n+k)p}{q}} = \frac{1}{\sqrt{B}}
$$
and
$$
\prod_{n \geq 0}
\prod_{\genfrac{}{}{0cm}{2}{0 < k < B}{k \not \equiv 0 \bmod q}}
\left(\frac{Bn+k}{Bn+k+1}\right)^{\sin \tfrac{\pi k p}{q}
\cos \tfrac{\pi(2n+k)p}{q}} = 1.
$$
\end{corollary}

\proof Let $\omega := e^{2 \pi ip/q}$. Since $B \equiv 1 \bmod q$, we may take
$u(n) := \omega^n$ in Addendum (2), yielding the formula
$$
\sum_{n \geq 0}
\sum_{\genfrac{}{}{0cm}{2}{0 < k < B}{k \not \equiv 0 \bmod q}}
\omega^n(1-\omega^k) \log \frac{Bn+k}{Bn+k+1} = - \log B.
$$
Writing
$\omega^n(1 - \omega^k) = -2i \omega^{n + \frac{k}{2}} \sin \frac{\pi k p}{q}$,
and multiplying the real and imaginary parts of the formula by $1/2$,
the result follows.  
\endpf

\begin{example}\label{iandsigma}
Take $B=5$, $p=1$, and $q=4$. Squaring the products, we get

{\it Define $\sigma(n)$ to be $+1$ if $n$ is a square modulo $4$, 
and $-1$
otherwise, that is,
$$
\sigma(n) :=
\begin{cases}
+1  &\mbox{\rm if } n \equiv 0 \mbox{\rm \ or } 1 \bmod 4, \\
-1  &\mbox{\rm if } n \equiv 2 \mbox{\rm \ or } 3 \bmod 4.
\end{cases}
$$
Then
$$
\prod_{n \geq 0} \left(\frac{5n+1}{5n+2}\right)^{\sigma(n)}
\left(\frac{5n+2}{5n+3}\right)^{\sigma(n)+\sigma(n+1)}
\left(\frac{5n+3}{5n+4}\right)^{\sigma(n+1)} = \frac{1}{5}
$$
and
$$
\prod_{n \geq 0} \left(\frac{5n+1}{5n+2}\right)^{\sigma(n-1)}
\left(\frac{5n+2}{5n+3}\right)^{\sigma(n-1)+\sigma(n)}
\left(\frac{5n+3}{5n+4}\right)^{\sigma(n)} = 1.
$$
}
\end{example}

\subsection{The sum-of-digits function $s_B(n)$}\label{section-sum}

Other products can also be obtained from Corollary~\ref{generalmult}. We give
three corollaries, each of which generalizes the Woods-Robbins formula 
$P = 1/\sqrt{2}$ in the Introduction. Recall that $s_B(n)$ denotes the sum 
of the $B$-ary digits of the integer $n$.

\begin{corollary}\label{zsB}
Fix an integer $B > 1$ and a complex number $z$ with $|z| \leq 1$.
If $z \not\in \{0,1\}$, then
$$ 
\prod_{n \geq 0} \prod_{\genfrac{}{}{0cm}{2}{0 < k < B}{z^k \neq 1}}
\left(\frac{Bn+k}{Bn+k+1}\right)^{z^{s_B(n)}(1-z^k)} = \frac{1}{B}\cdot
$$
\end{corollary}

\proof Take $u(n) := z^{s_B(n)}$ in Corollary~\ref{generalmult}
and note that $s_B(k) = k$ when $0 < k < B$. \endpf

\begin{example}
Take $B=2$ and $z=1/2$. Squaring the product, we obtain
$$
\prod_{n \geq 0} \left(\frac{2n+1}{2n+2}\right)^{(1/2)^{s_2(n)}} = \frac{1}{4}\cdot
$$
\end{example}

\begin{corollary}\label{rootandsum}
Let $B, p, q$ be integers with $B > 1$ and $q > p > 0$. Then
$$
\prod_{n \geq 0}
\prod_{\genfrac{}{}{0cm}{2}{0 < k < B}{k \not \equiv 0 \bmod q}}
\left(\frac{Bn+k}{Bn+k+1}\right)^{\sin \tfrac{\pi k p}{q} 
\sin \tfrac{\pi(2s_B(n)+k)p}{q}} = \frac{1}{\sqrt{B}}
$$
and
$$
\prod_{n \geq 0}
\prod_{\genfrac{}{}{0cm}{2}{0 < k < B}{k \not \equiv 0 \bmod q}}
\left(\frac{Bn+k}{Bn+k+1}\right)^{\sin \tfrac{\pi k p}{q}
\cos \tfrac{\pi(2s_B(n)+k)p}{q}} = 1.
$$
\end{corollary}

\proof Use the proof of Corollary~\ref{rootofunity}, but replace $B \equiv 1 \bmod q$
with $s_B(Bn+k) = s_B(n)+k$ when $0 \leq k < B$, and replace $\omega^n$ with 
$\omega^{s_B(n)}$.  \endpf

\begin{example}
Take $B=2$, $q=4$, and $p=1$. Squaring the products and defining 
$\sigma(n)$ as in Example~\ref{iandsigma}, we get
$$
\prod_{n \geq 0} \left(\frac{2n+1}{2n+2}\right)^{\sigma(s_2(n))} 
= \frac{1}{2} \ \ \ \mbox{and} \ \ \
\prod_{n \geq 0} \left(\frac{2n+1}{2n+2}\right)^{\sigma(s_2(n)+1)} = 1.
$$
\end{example}

\bigskip

In the same spirit, we recover a result from \cite[p. 369-370]{AMFP}. 

\begin{example}\label{theta}

Taking $B=q=3$ and $p=1$ in Corollary~\ref{rootandsum}, we obtain two 
infinite products. Raising the second to the power $-2/\sqrt{3}$ and 
multiplying by the square of the first, we get

{\em Define $\theta(n)$ by
$$
\theta(n) :=
\begin{cases}
\ \ 1 &\mbox{\rm if } n \equiv 0 \ \mbox{\rm or } 1 \bmod 3, \\
-2    &\mbox{\rm if } n \equiv 2 \bmod 3.\\
\end{cases}
$$
Then
$$
\prod_{n \geq 0} 
(3n+1)^{\theta(s_3(n))} (3n+2)^{\theta(s_3(n)+1)} (3n+3)^{\theta(s_3(n)+2)} 
= \frac{1}{3}.
$$
}
\end{example}

\begin{corollary}[\cite{Sondow}]\label{sumofdigits}
Let $B > 1$ be an integer. Then
$$
\prod_{n \geq 0}
\prod_{\genfrac{}{}{0cm}{2}{0 < k < B}{k \ \mbox{\scriptsize{\rm odd}}}}
\left(\frac{Bn+k}{Bn+k+1}\right)^{(-1)^{s_B(n)}}
= \frac{1}{\sqrt{B}}\cdot
$$
\end{corollary}

\proof Take $z=-1$ in Corollary~\ref{zsB} 
(or take $q = 2$ and $p=1$ in Corollary~\ref{rootandsum}).  \endpf

\begin{example}
With $B = 2$, since $s_2(n) = N_{1,2}(n)$, we recover the Woods-Robbins 
formula $P = 1/\sqrt{2}$. Taking $B=6$ gives
$$
\prod_{n \geq 0} \left(
\frac{(6n+1)(6n+3)(6n+5)}{(6n+2)(6n+4)(6n+6)}
\right)^{(-1)^{s_6(n)}} = \frac{1}{\sqrt 6}\cdot
$$
\end{example}

\begin{remark}\label{trickson}
Corollary~\ref{sumofdigits} can also be obtained from \cite[Theorem~1]{ACMS},
as follows. Taking $x$ equal to $-1$ and $j$ equal to $0$ in that theorem gives
$$
\sum_{n \geq 0} (-1)^{s_B(n)} \log \frac{n+1}{B \lfloor n/B \rfloor +B} = 
-\frac{1}{2} \log B
$$
where $\lfloor x \rfloor$ is the integer part of $x$.
But the series is equal to
$$
\begin{array}{lll}
\displaystyle\sum_{m \geq 0} \sum_{0 \leq k < B}
(-1)^{s_B(Bm+k)} \log \frac{Bm+k+1}{Bm+B}
&=& \displaystyle\sum_{m \geq 0} (-1)^{s_B(m)} \sum_{0 \leq k < B} (-1)^k
\log \frac{Bm+k+1}{Bm+B} \\
&=& \displaystyle\sum_{m \geq 0} (-1)^{s_B(m)}
\sum_{\genfrac{}{}{0cm}{2}{k \ \mbox{\scriptsize{\rm odd}}}{0 < k < B}}
\log \frac{Bm+k}{Bm+k+1}
\end{array}
$$
where the last equality follows by looking separately at the cases $B$ even
and $B$ odd.
\end{remark}

\subsection{The counting function $N_{j,B}(n)$}

We can also compute some infinite products associated with counting the number
of occurrences of one or several given digits in the base $B$ expansion of 
an integer.

\begin{definition}
If $B$ is an integer $\geq 2$ and if $j$ is in $\{0, 1, \ldots, B-1\}$, let 
$N_{j,B}(n)$ be the number of occurrences of the digit $j$ in the $B$-ary 
expansion of $n$ when $n > 0$, and set $N_{j,B}(0) := 0$. 
\end{definition}

\begin{corollary}\label{suppl}
Let $B, q, p$ be integers with $B > 1$ and $q > p > 0$. 
Let $J$ be a nonempty, proper subset of $\{0, 1, \ldots, B-1\}$. Define 
$N_{J,B}(n) := \sum_{j \in J} N_{j,B}(n)$. Then the following equalities hold: 
$$
\prod_{k \in J}
\prod_{n \geq \delta_k} 
\left(\frac{Bn+k}{Bn+k+1}\right)^{\sin\tfrac{\pi(2N_{J,B}(n)+1)p}{q}} 
= B^{-\tfrac{1}{2\sin\tfrac{\pi p}{q}}}
$$
and
$$
\prod_{k \in J}
\prod_{n \geq \delta_k} 
\left(\frac{Bn+k}{Bn+k+1}\right)^{\cos\tfrac{\pi(2N_{J,B}(n)+1)p}{q}}
= 1.
$$
\end{corollary}

\proof  Let $\omega := e^{2 \pi i p/q}$. We denote 
$u_{q,j,B}(n) := \omega^{N_{j,B}(n)}$ 
and $u_{q,J,B}(n) := \prod_{j \in J} u_{q,j,B}(n) = \omega^{N_{J,B}(n)}$.
Note that, for every $j$ in $\{1, 2, \ldots, B-1\}$, the sequence
$(u_{q,j,B}(n))_{n \geq 0}$ is strongly $B$-multipli\-ca\-tive and 
nonzero, hence satisfies Hypothesis ${\cal H}_B$. The sequence 
$(u_{q,0,B}(n))_{n \geq 0}$ also satisfies Hypothesis ${\cal H}_B$, 
as is seen by taking $z = \omega$ in Example~\ref{special}. 
Therefore the sequence $(u_{q,J,B}(n))_{n \geq 0}$ satisfies Hypothesis 
${\cal H}_B$, with, for $k = 0, 1, \ldots, B-1$, the value $v(k) := \omega$ 
if $k \in J$ and $v(k) := 1$ otherwise. 

Now $|u_{q,J,B}(n)| = 1$ for $n \geq 0$, and $|v(k)| = 1$ for 
$k = 0, 1, \ldots, B-1$. Furthermore, $|\sum_{0 \leq k < B} v(k) | < B$, 
since $v$ is not constant on $\{0, 1, \ldots, B-1\}$.
Thus we may apply Addendum (1) with $u(n) := u_{q,J,B}(n)$, obtaining 
$$
(1 - \omega) \sum_{k \in J} \sum_{n \geq \delta_k} 
\omega^{N_{J,B}(n)} \log\frac{Bn+k}{Bn+k+1} = - \log B.
$$
Writing $(1-\omega)\omega^{N_{J,B}(n)} =
-2i\omega^{N_{J,B}(n)+\frac{1}{2}}\sin \frac{\pi p}{q}$, and taking the 
real and imaginary parts of the summation, the result follows.  \endpf

\begin{example}\label{J-set}
Taking $q=2$ and $p=1$ in the first formula gives
$$
\prod_{k \in J}
\prod_{n \geq \delta_k}
\left(\frac{Bn+k}{Bn+k+1}\right)^{(-1)^{N_{J,B}(n)}} = \frac{1}{\sqrt{B}}\cdot
$$
An application is an alternate proof of Corollary~\ref{sumofdigits}: 
take $J$ to be the set of odd numbers in $\{1, 2, \ldots, B-1\}$; 
since $s_B(n) = \sum_{0<k<B} kN_{k,B}(n)$, it follows that 
$(-1)^{\sum_{j \in J} N_{j,B}(n)} = (-1)^{s_B(n)}$.
\end{example}

\begin{remark}
Corollary~\ref{suppl} requires that $J$ be a proper subset of
$\{0, 1, \ldots, B-1\}$. Suppose instead that  $J = \{0, 1, \ldots, B-1\}$. 
Then $N_{J,B}(n)$ is the number of $B$-ary digits of $n$ if $n > 0$ 
(that is, $N_{J,B}(n) = \lfloor \frac{\log n}{\log B} \rfloor + 1$), and 
$N_{J,B}(0) = 0$.  In that case, Corollary~\ref{suppl} does not apply, and 
the products may diverge. For example, when $B=q=2$ and $p=1$ the logarithm 
of the first product is equal to the series 
$$
- \log 2 + \sum_{n \geq 1} (-1)^{\left\lfloor\tfrac{\log n}{\log 2}\right\rfloor}
\log \frac{n+1}{n},
$$
which does not converge. However, note its resemblance with
Vacca's (convergent) series for Euler's constant \cite{Vacca}
$$
\gamma = \sum_{n \geq 1} \left\lfloor\frac{\log n}{\log2}\right\rfloor 
\frac{(-1)^n}{n}.
$$
\end{remark}

\begin{corollary}\label{singledigit}
Let $B, q, p$ be integers with $B > 1$ and $q > p > 0$.
Then for $k = 0, 1, \ldots, B-1$ the following equalities hold:
$$
\prod_{n \geq \delta_k}
\left(\frac{Bn+k}{Bn+k+1}\right)^{\sin\tfrac{\pi(2N_{k,B}(n)+1)p}{q}}
= B^{-\tfrac{1}{2\sin\tfrac{\pi p}{q}}}
$$
and
$$
\prod_{n \geq \delta_k}
\left(\frac{Bn+k}{Bn+k+1}\right)^{\cos\tfrac{\pi(2N_{k,B}(n)+1)p}{q}}
= 1.
$$
\end{corollary}

\proof Take $J := \{k\}$ in Corollary~\ref{suppl}. (The case $k=0$ and 
$p=1$ is Example~\ref{counting0} with $z = e^{2\pi i/q}$.)  \endpf

\begin{example}\label{q2p1} 
Taking $q=2$ and $p=1$ in the first formula (or taking $J = \{k\}$ in 
Example~\ref{J-set}) yields
$$
\prod_{n \geq \delta_k}
\left(\frac{Bn+k}{Bn+k+1}\right)^{(-1)^{N_{k,B}(n)}} = \frac{1}{\sqrt{B}}\cdot
$$
In particular, if $B=2$ the choice $k=1$ gives the Woods-Robbins formula
$P = 1/\sqrt{2}$, and $k=0$ gives 
$$
\prod_{n \geq 1} \left(\frac{2n}{2n+1}\right)^{(-1)^{N_{0,2}(n)}} 
= \frac{1}{\sqrt{2}}\cdot
$$
\end{example}

\begin{remark}
For base $B=2$, the formulas in Example~\ref{q2p1} are special cases of 
results in \cite{AS89}, where $N_{j,2}(n)$ is generalized to counting 
the number of occurrences of a given {\it word} in the binary expansion 
of $n$. On the other hand, the value of the product $Q$ in the Introduction,
$$
Q = \prod_{n \geq 1} \left(\frac{2n}{2n+1}\right)^{(-1)^{N_{1,2}(n)}},
$$
remains a mystery.
\end{remark}

\begin{example}
Take $B=q=3$ and $p=1$. Raising the first product to the power 
$2/\sqrt{3}$ and squaring the second, we obtain

{\em Define $\eta(n)$ by
$$
\eta(n) :=
\begin{cases}
+1 &\mbox{\rm if } n \equiv 0  \bmod 3, \\
\ \ 0 &\mbox{\rm if } n \equiv 1 \bmod 3,\\
-1    &\mbox{\rm if } n \equiv 2 \bmod 3,\\
\end{cases}
$$
and define $\theta(n)$ as in Example~\ref{theta}.
Then for $k = 0, 1$, and $2$
$$
\prod_{n \geq \delta_k} \left(\frac{3n+k}{3n+k+1}\right)^{\eta(N_{k,3}(n))} 
= \frac{1}{3^{2/3}}
\ \ \ \mbox{\it and} \ \ \
\prod_{n \geq \delta_k} \left(\frac{3n+k}{3n+k+1}\right)^{\theta(N_{k,3}(n)+1)}
= 1.
$$
}
\end{example}

\subsection{The Gamma function}

It can happen that the exponent in some of our products is a periodic 
function of $n$. For example, this is obviously the case in 
Corollary~\ref{rootofunity}. To take another example, it is not hard to 
see that if $B$ odd, then $(-1)^{s_B(n)} = (-1)^n$. Hence 
Corollary~\ref{sumofdigits} gives
\renewcommand{\theequation}{\fnsymbol{equation}}
\setcounter{equation}{6}
\begin{equation}\label{Bodd}
\prod_{n \geq 0}
\prod_{\genfrac{}{}{0cm}{2}{0 < k <B}{k \ \mbox{\scriptsize{\rm odd}}}}
\left(\frac{Bn+k}{Bn+k+1}\right)^{(-1)^n}
= \frac{1}{\sqrt{B}} \ \ \mbox{\rm ($B$ odd)}.
\end{equation}
(This formula can also be obtained from Corollary~\ref{rootofunity} 
with $q = 2$ and $p=1$.) For instance
$$
P_{1,3} := \prod_{n \geq 0} \left(\frac{3n+1}{3n+2}\right)^{(-1)^n} 
= \frac{1}{\sqrt{3}}.
$$

The product $P_{1,3}$ can also be computed using the following 
corollary of the Weierstrass product for the Gamma function
\cite[Section~12.13]{WW}.

{\it If $d$ is a positive integer and $a_1 + a_2 + \cdots + a_d = 
b_1 + b_2 + \cdots + b_d$, where the $a_j$ and $b_j$ are complex numbers
and no $b_j$ is zero or a negative integer, then 
$$
\prod_{n \geq 0}\frac{(n+a_1) \cdots (n+a_d)}{(n+b_1) \cdots (n+b_d)}
= \frac{\Gamma(b_1) \cdots \Gamma(b_d)}{\Gamma(a_1) \cdots \Gamma(a_d)}.
$$
}
Combining this with the relation $\Gamma(x) \Gamma(1-x) = \pi/\sin \pi x$
gives $P_{1,3} = 1/\sqrt{3}$. 

The computation can be generalized, using Gauss' multiplication theorem 
for the Gamma function, to give another proof of Corollary~\ref{sumofdigits} 
for $B$ odd. Likewise, an analog of the odd-$B$ case of 
Corollary~\ref{sumofdigits} can be proved for even $k$:
$$
\prod_{n \geq 1}
\prod_{\genfrac{}{}{0cm}{2}{0 \leq k < B}{k \ \mbox{\scriptsize{\rm even}}}}
\left(\frac{Bn+k}{Bn+k+1}\right)^{(-1)^n}
= \frac{\pi \sqrt{B}}{2^B} \binom{B-1}{(B-1)/2} \ \ \mbox{\rm ($B$ odd)}.
$$
Multiplying this by the formula
$$
\prod_{n \geq 1}
\prod_{\genfrac{}{}{0cm}{2}{0 < k < B}{k \ \mbox{\scriptsize{\rm odd}}}}
\left(\frac{Bn+k}{Bn+k+1}\right)^{(-1)^n}
= \frac{2^{B-1}}{\sqrt{B}} \binom{B-1}{(B-1)/2}^{-1} \ \ \mbox{\rm ($B$ odd)},
$$
which is (\ref{Bodd}) rewritten, yields Wallis' product for $\pi$. (For an 
evaluation of the preceding two products when $B = 2$, see \cite[Example~7]{SH}.)

\end{document}